\RequirePackage{fix-cm}
\documentclass{svjour3}                     
\smartqed  
\usepackage{url}
\usepackage{mathptmx}      
\usepackage{amssymb}
\usepackage{amsmath}
\usepackage{wasysym}
\usepackage{graphicx}
\usepackage[font=footnotesize]{subcaption}
\usepackage{animate}

\renewcommand{\emph}[1]{\textsl{#1}}
\newcommand{\bi}{\begin{itemize}}
\newcommand{\ei}{\end{itemize}}
\newcommand{\bc}{\begin{center}}
\newcommand{\ec}{\end{center}}
\newcommand{\bea}{\setlength\arraycolsep{2pt}\begin{eqnarray}}
\newcommand{\eea}{\end{eqnarray}\setlength\arraycolsep{6pt}}
\newcommand{\s}{\section}
\newcommand{\subs}{\subsection}
\newcommand{\ie}{\emph{i.e.}}

\newcommand{\M}{\textsc{Matlab}}

\DeclareMathOperator*{\argmin}{argmin}

\journalname{Numerical Algorithms}

\begin{document}

\title{Optimal Residuals and the Dahlquist Test Problem}
\author{Robert~M.~Corless${}^1$, C.~Yal\c{c}{\i}n~Kaya${}^2$ and Robert~H.~C.~Moir${}^1$}
\institute{
	${}^1$The Rotman Institute of Philosophy and\\
         The Ontario Research Center for Computer Algebra and\\
         The School of Mathematical and Statistical Sciences\\
	The University of Western Ontario\\
	${}^2$School of Information Technology and Mathematical Sciences\\
	University of South Australia
}

\date{Received: date / Accepted: date}




\maketitle
\vspace{-60pt}
\bc
\includegraphics[scale=0.15]{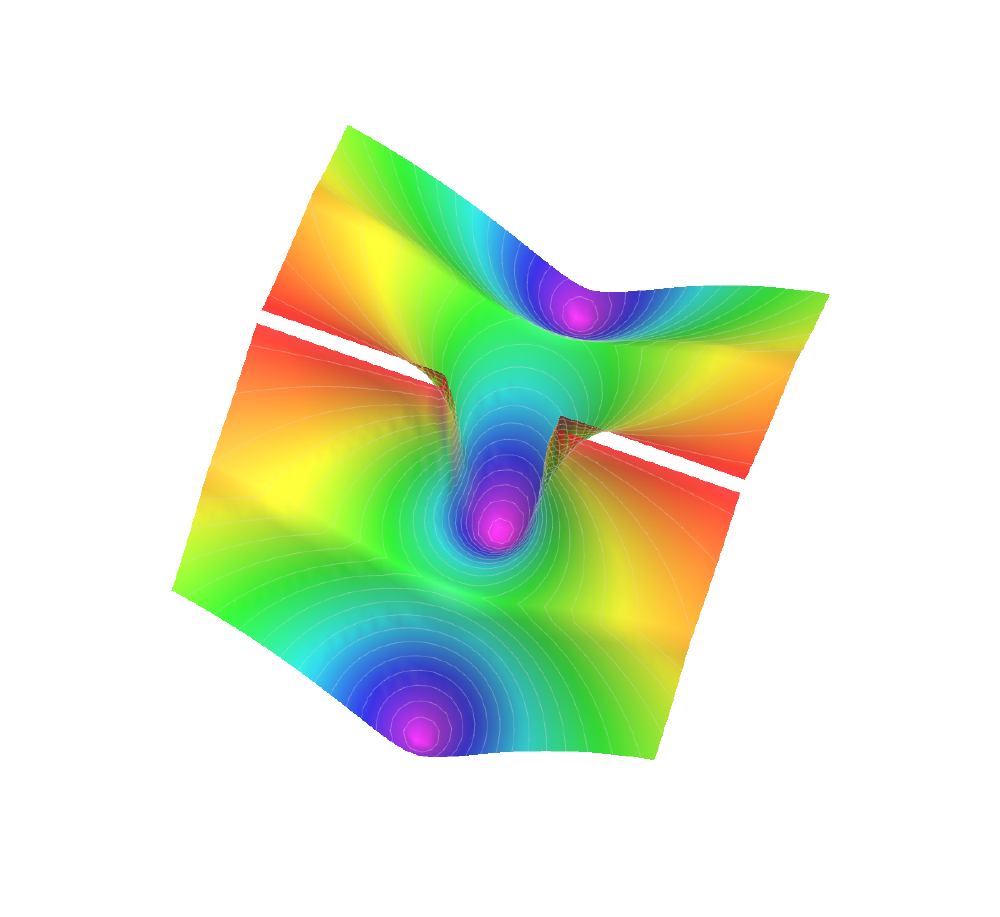}
\ec

\begin{abstract}
We show how to compute the \emph{optimal relative backward error} for the numerical solution of the Dahlquist test problem by one-step methods.
This is an example of a general approach that uses results from optimal control theory to compute optimal residuals, but elementary methods can also be used here because the problem is so simple.  This analysis produces new insight into the numerical solution of stiff problems. 
\end{abstract}

\s{Introduction}
The study of stiff differential equations and their efficient numerical solution is by now a mature field. There are several, perhaps many, efficient practical methods with freely available high quality implementations. The literature on the theory of such methods is extensive. Surprisingly, it is not yet complete: for instance, see the survey \cite{soderlind2014stiffness}, which has the intriguing title \emph{Stiffness 1952-2012: Sixty years in search of a definition}. That paper re-examines the fundamentals and thoroughly surveys the literature, and proposes a new stiffness indicator that they claim is useful both \emph{a priori} for indicating stiffness and \emph{a posteriori} for indicating varying regions of the solution that stiffness is important.

This paper takes a different approach, that of \emph{optimal residuals}, \emph{i.e.,}\ backward error, and uses it on the Dahlquist test problem to generate some new observations about this, the simplest of all stiff problems. Indeed, \cite{soderlind2014stiffness} calls this problem ``simplistic'' and with good reason, but surprisingly it still has things to teach us.

Trying to study stiff problems from the point of view of backward error analysis is itself not new. For instance, there is the PhD thesis of W.L. Seward and the paper~\cite{enright1989achieving}. But there is an intrinsic dissonance: a stiff problem has the feature that errors are (often sharply) damped as the integration proceeds forward in time, and thus it is not obvious why one might prefer to look at backward error when, if the right method is used, small forward error happens more or less automatically.

We claim that backward error, in particular optimal residuals, which we explain below, really can be useful. One way to see the usefulness of this type of analysis is in contrast to the classical stability analysis on the Dahlquist test problem, $\dot{y}=\lambda y$, $\lambda\in\mathbb{C}$. In the classical stability analysis there is an emphasis on A-stability, with an understanding that, at least locally where linear stability assessments are valid, a basic criterion for stability of a method is that the numerical solution decays where the actual solution does. This condition is of particular concern for stiff problems, since it entails that the method will not be subject to stability restrictions on account of eigenvalues with large negative real part.

The classical stability analysis really does explain most of the behaviour of implicit methods on stiff problems, and much insight has been gained thereby.  This present analysis is a refinement only, that offers the possibility of explaining some extra considerations which are ``well-known'' to practitioners, namely that, under certain circumstances, a stiff method may well be stable but not as accurate as one might wish.
With the optimal residual approach a rather different set of criteria emerge, which imply a de-emphasis of the criterion of A-stability, being enhanced with a consideration of those regions of the problem space that a given method can reproduce the dynamics of the reference problem accurately; such regions often extend comfortably into the right-half of the complex plane, even for explicit methods. 



There are also deeper reasons why backward error analysis is useful on stiff problems, as a result of the fact that a small backward error entails a small perturbation of the dynamics of the problem. In this regard, one place where optimal residuals are especially useful is in the solution of systems which have nontrivial attracting sets: the problems can be stiff because decay to the attracting set can be very strong, but small backward error is important too in order to get the dynamics right on the attracting set, which might in fact be actually chaotic.  Of course, for a chaotic problem, small forward error is not possible at all (at least, not for very long), but as explained for instance in~\cite{corless1994good} small backward error can be perfectly satisfactory as an explanation of the success of a numerical method on a chaotic problem.

The question of what good a numerical method is, even one providing a small backward error, for computation of an attracting set for a chaotic system, deserves a repetition of the answer given in that reference. Every numerical analyst knows that small forward error requires exponentially accurate initial conditions and exponentially accurate integration. One way to explain a successful computation is if some form of shadowing is invoked. That is, the forward accuracy of the computed trajectory is explained by being``shadowed'' by an exact solution of the reference problem, typically from some nearby initial condition.  This is a form of backward error analysis~\cite{moir2010reconsidering}.  Computationally verifying that shadowing has occurred is expensive, however, and while shadowing is generic, there are no \emph{a priori} guarantees that shadowing will occur or has occurred.

A cheaper and simpler method to explain the success of a numerical method on a chaotic problem is to verify that the \emph{residual} (also called the defect, deviation, or slope error) is small. 
Also take note of the utility of interpreting ordinary numerical errors as a modelling error in many circumstances. If the model equations are written
$$\dot{y}=f(t,y)$$
after having been derived for some physical context, universally by neglecting minor effects, and our numerical solution gives us $z(t)$ exactly satisfying
$$\dot{z}(t)=f(t,z(t))+\Delta(t).$$
or perhaps
$$\dot{z}(t)=(I+\delta(t))f(t,z(t)),$$
and $\Delta$ or $\delta$ are small compared to the neglected terms, then as Nick Higham puts it ``for all we know, we already have the exact solution,'' for error has been introduced into the dynamics by modelling assumptions.

We also emphasize that there \emph{must} be some feature of this system that is robust under perturbations, or else even the reference solution of $\dot{y}=f(t,y)$ would be useless in the face of real-life perturbations. The existence of such a feature for a given system was termed ``well-enough conditioning'' in~\cite{corless1994error}. Nearly all models used have this property, even chaotic ones. For instance, the dimension of the Lorenz attractor is quite robust under forcings of this type: the attracting sets of 
\bea
\dot{y_1}&=&\sigma(y_2-y_1)+\delta_1(y),\nonumber\\ 
\dot{y_2}&=&y_1(\rho-y_3)-y_2+\delta_2(y),\nonumber\\
\dot{y_3}&=&y_1y_2-\beta y_3+\delta_3(y),\nonumber
\eea
are remarkably close to one another even for quite large $\Delta(t)=(\delta_1(t),\delta_2(t),\delta_3(t))^{\text{T}}$.~\cite{corless1994error}

So backward error may be important in explaining the success or failure of numerical methods for chaotic systems, which can be stiff. More than this, however, we believe our approach to stability analysis on the Dahlquist test problem may lead to a refinement of the classical explanation of the practical success of various numerical methods. 

We now consider backward error and the general context of\emph{optimal residuals} briefly, before moving on to an elementary analysis of the Dahlquist test problem.

\s{Optimal Residuals and Backward Error}

%
%
%
%

\subs{Backward Error, Residuals and Interpolants}
Most codes supply interpolants: for graphical output, for event location, or for handling delays. Given an interpolant, which we will call $z(t)$, to the numerical \textsl{skeleton} $(t_k,y_k)$ of our computed solution to the ODE $\dot y = f(t,y)$, we define the \textsl{residual} $\Delta(t)$ as
\begin{equation}
\Delta(t) := \dot z(t) - f( t, z(t)) \>. \label{eq:residualdef}
\end{equation}
As previously  stated, this is sometimes called the defect. In one sense, this is a kind of backward error, because the computed $z(t)$ can be interpreted as the exact solution of the perturbed equation
\begin{equation}
\dot y(t) = f(t, y(t) ) + \Delta(t)\>. \label{eq:perturbedeq}
\end{equation}
Note that the residual as defined here is dependent on just which interpolant $z(t)$ is used.

To be compatible with the accuracy of the numerical solutions they interpolate, these interpolants should be $O(h^p)$ accurate, but \emph{sometimes are not}. For example,  in \M, \texttt{ode45} uses a fifth-order Runge-Kutta Fehlberg formula, but has only a fourth-order interpolant: so $z'(t)$ will only be third-order accurate. A consequence is that the residual, or defect, $\Delta(t):=\dot{z}(t)-f(t,z)$ will sometimes be overestimated.

The work of \cite{enright1989new,enright2000continuous}, \cite{enright1991parallel}, \cite{enright2007robust}, \cite{muir2003pmirkdc}, \cite{cao2004posteriori} has shown that one can interpolate the skeleton $\{t_k,y_k\}$ of a numerical solution in a practical way that gives the correct asymptotic size of the residual $\|\Delta\|=O(h^p)$ as the mean time step $h\rightarrow 0$, for a method of order $p$. This gives a practical method that gives tolerance proportionality and robust reliability for tight tolerances.


If instead one interpolates the skeleton badly, one will obviously get a large residual.  In one of their examples, Hubbard \& West wonder if a smaller residual (they term it ``slope error'') might be achieved with some other interpolant~\cite{HubbardWest:1991}. This raises the question of how to find the ``best'' interpolant, which gives the ``smallest'' residual.  This smallest residual is the one that will most accurately  measure how good a job the underlying method did in generating the skeleton.

This question is answered in general elsewhere (Corless, Kaya \& Moir, in preparation) where optimization methods are used to find interpolants from $(t_k,y_k)$ to $(t_{k+1},y_{k+1})$ that minimize $\|\Delta(t)\|$ or $\|\delta(t)\|$ (either the $2$-norm or $\infty$-norm are handled). 

\subsection{Interpolants via optimal control theory}

For ease of exposition in this section we work entirely in $\mathbb{R}^n$.

Define the {\em control variable vector} $u(t) := \delta(t)$, where $u:[t_i,t_{i+1}]\to\mathbb{R}^n$, as a piecewise continuous vector function. The problem of finding interpolants $z:[t_i,t_{i+1}]\to\mathbb{R}^n$ from $(t_i,y_i)$ to $(t_{i+1},y_{i+1})$ that minimize the $L^\infty$-norm of the relative error $\|\delta\|_{L^\infty}$ can be stated as an optimal control problem:
\[
\mbox{P1: }\left\{\begin{array}{rl}
\displaystyle\mbox{minimize} & \ \ \displaystyle\,\max_{t_i\le t\le t_{i+1}} \|u(t)\|_{\infty} \\[3mm]
\mbox{subject to} & \ \ \dot{z}(t) = f(t,z(t))\,(1 + u(t))\,, \mbox{ a.e. } t_i\le t\le t_{i+1}\,, \\[2mm]
& \ \ z(t_i) = y_i,\ \ z(t_{i+1}) = y_{i+1}\,,
\end{array} \right.
\]
where $\|\cdot\|_\infty$ is the $\ell_\infty$-norm in $\mathbb{R}^n$.  Problem~(P1), where $z(t)$ assumes the role of a {\em state variable vector},
can be transformed into a more standard form to apply a maximum principle, as in \cite{KayNoa2013}, as follows. Let $\alpha$ be a new parameter.  Then one can rewrite Problem~(P1) equivalently as
\[
\mbox{P2: }\left\{\begin{array}{rl}
\displaystyle\mbox{minimize} & \ \ \alpha \\[3mm]
\mbox{subject to} & \ \ \dot{z}(t) = f(t,z(t))\,(1 + u(t))\,, \mbox{ a.e. } t_i\le t\le t_{i+1}\,, \\[2mm]
& \ \ z(t_i) = y_i,\ \ z(t_{i+1}) = y_{i+1}\,, \\[2mm]
& \ \ |u_j(t)| \le \alpha\,,\ \mbox{ a.e. } t_i\le t\le t_{i+1}\,,\ j = 1,\ldots,n\,.
\end{array} \right.
\]
To apply the maximum principle, we define the {\em Hamiltonian function} for Problem~(P2):
\[
H(z,u,\psi,t) := \sum_{j=1}^n\,\psi_j\,f_j(t,z)\,(1 + u_j)\,,
\]
where $\psi:[t_i,t_{i+1}]\to\mathbb{R}^n$ is referred to as the {\em adjoint} (or {\em costate}) {\em variable vector}, such that
\begin{equation}  \label{adjoint}
\dot{\psi}_j = -\frac{\partial H}{\partial z_j} = -\sum_{j=1}^n\,\psi_j\,\frac{\partial f_j}{\partial z_j}(t,z)\,(1 + u_j)\,,\ j = 1,\ldots,n\,.
\end{equation}
The maximum principle asserts that if $u$ is optimal then there exists a nontrivial adjoint variable vector such that $u$ minimizes the Hamiltonian; namely
\[
u(t)\in\argmin_{|v|\le \alpha} H(z(t),v,\psi(t),t)\,,\ \mbox{a.e. } t\in[t_i,t_{i+1}]\,,
\]
i.e.,
\begin{equation}  \label{control}
u_j(t) = \left\{\begin{array}{ll}
\ \ \alpha\,, & \mbox{if}\  \psi_j(t)\,f_j(t,z(t)) < 0\,, \\[1mm]
-\alpha\,, & \mbox{if}\ \psi_j(t)\,f_j(t,z(t)) > 0\,, \\[1mm]
\mbox{\small undetermined}\,, & \mbox{if}\ \psi_j(t)\,f_j(t,z(t)) = 0\,.
\end{array}\right.\quad j = 1,\ldots,n\,.
\end{equation}
This implies that the components of the optimal control, or the optimal backward error, will either be {\em bang--bang} (i.e., switching between the constants $\alpha$ and $-\alpha$) or {\em singular} (when $\psi_j(t)\,f_j(t,z(t)) = 0$ over some nontrivial subinterval).  It is not uncommon to encounter the error components to be of bang--bang or singular types as preliminary studies show.  However, this is a topic of further research.

In this paper, we focus on the {\em Dahlquist test equation}, $\dot{y}=\lambda y$, for which the optimality conditions derived above simplify drastically. For Dahlquist with real $\lambda$, the adjoint equation~\eqref{adjoint} reduces to the scalar ODE
\begin{equation}  \label{adjoint2}
\dot{\psi}(t) = -\lambda\,(1 + u(t))\,\psi(t)\,.
\end{equation}
Singular control can be ruled out for this case as follows.  If the optimal control is singular, then from \eqref{control} either $\psi(t) = 0$ or $f(t,z(t)) = z(t) = 0$ over a subinterval $[s_1,s_2]\subset[t_i,t_{i+1}]$. Suppose that $\psi(t) = 0$ over $[s_1,s_2]$.  Then $\dot{\psi}(t) = -\lambda\,(1 + u(t))\,\psi(t) = 0$ over $[s_1,s_2]$, for any $u$.  However, this further means that $\psi(t) = 0$ for all $t\in[t_i,t_{i+1}]$, which, being trivial, is not allowed by the maximum principle. Now, suppose that $z(t) = 0$ over $[s_1,s_2]$.  Then similarly $z(t) = 0$ for all $t\in[t_i,t_{i+1}]$, which we can leave out  as a very special case with $y_{i} = 0$.  

Therefore the only possible optimal relative error is of bang--bang type. In other words, $|\delta(t)| = \alpha$, i.e., the optimal relative error is of constant magnitude.

Suppose that the control $u(t) = \alpha$ or $-\alpha$ on $[t_i,t_i+\varepsilon]$, with a small $\varepsilon>0$.  Then \eqref{adjoint2} implies that $\psi(t)$ does not change sign for any $t\in[t_i,t_{i+1}]$.  Thus, either $u(t) = \alpha$ or $u(t) = -\alpha$ for all $t\in[t_i,t_{i+1}]$. As a result, the optimal interpolant is given by
\[
\dot{z}(t) = (1\pm\alpha)\,\lambda\,z(t)\,,
\]
whose solution can simply be obtained as
\[
z(t) = y_i\,e^{(1\pm\alpha)\,\lambda\,(t-t_i)}\,.
\]
It follows that
\begin{equation}
y_{i+1} = y_i\,e^{(1\pm\alpha)\,\lambda\,h}\,, \label{eq:expon}
\end{equation}
or\footnote{Going from equation~\eqref{eq:expon} to equation~\eqref{eq:logarithm} is correct over the positive reals but not \emph{quite} right over $\mathbb{C}$: $\ln z$ means the principal branch of the logarithm, with argument in $(-\pi,\pi]$. We will see how to choose the correct complex branch in equation~\eqref{eq:firstunwind}, in a later section.}
\begin{equation}
1\pm\alpha = \frac{1}{\lambda\,h}\,\ln\frac{y_{i+1}}{y_i}\,.
\label{eq:logarithm}
\end{equation}
Rearranging further yields the minimum value of the $L^\infty$-norm of the relative error as
\[
\alpha = \left|1 - \frac{1}{\lambda\,h}\,\ln\frac{y_{i+1}}{y_i} \right| .
\]

One may also look for a control defined over more than one piece of the skeleton; this would 
allow the analysis to be used for multistep methods.  In this paper, however, we restrict attention to one-step methods.

Indeed for some problems the minimal $\delta$ is surprisingly easy to find: see also~\cite{corless2016variations} for an investigation of Torricelli's law using this idea.

We will use this approach to re-examine the Dahlquist test problem, $\dot{y}=\lambda y$. We will be able to find the interpolant giving the optimal $\delta$ without difficulty using only elementary methods. Our analysis will be valid for any one-step method. The results are quite surprising.

\s{Optimal Backward Error for $\dot{y}=\lambda y$}
Without loss of generality suppose that our one-step numerical method gives, for $\dot y  = \lambda y$, $y(t_n)=y_n$, and time-step $h$
$$y_{n+1}=y_n+h\Phi(y_n)=R(\mu)y_n\>,$$
where $\mu=\lambda h \in \mathbb{C}$ and the time-step $h > 0$. For instance, the Explicit Euler method gives
$$y_{n+1}=y_n+h\cdot\lambda y_n=(1+\mu)y_n,$$
while Implicit Euler gives
$$y_{n+1}=y_n+h\cdot\lambda y_{n+1}$$
or
$$y_{n+1}=\frac{1}{1-\mu}y_n.$$

A small table of $R(\mu)$ corresponding to various methods is given in Table \ref{tab:exp-approx}. Note that Taylor series methods of order $p$ have
$$R(\mu)=\sum_{k=0}^p\frac{\mu^k}{k!}.$$

\subsection{Classical stability regions}
The classical stability regions are defined given the behaviour of numerical methods on the Dahlquist test problem. The \emph{absolute stability region} is defined as the region in the complex plane where the value of $\mu=\lambda h$ guarantees that the Dahlquist test problem is uniformly bounded for all forward time steps. Again the rationale for focusing on such a simple problem is that it characterizes the local behaviour of a method on a scalar component of a linearized version of any nonlinear problem.

\begin{figure}[th!]
	\centering
	\begin{subfigure}[b]{0.3\textwidth}
		\centering
		\includegraphics[scale=0.2]{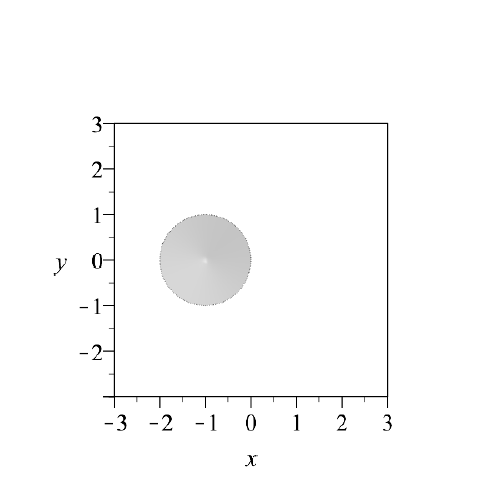}
		\caption{Euler method}
		\label{fig:theta-methods:a}
	\end{subfigure}
	\begin{subfigure}[b]{0.3\textwidth}
		\centering
		\includegraphics[scale=0.2]{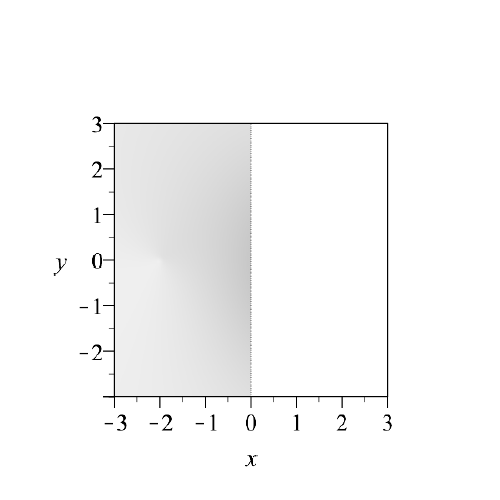}
		\caption{Implicit midpoint rule}
		\label{fig:theta-methods:b}
	\end{subfigure}
	\begin{subfigure}[b]{0.3\textwidth}
		\centering
		\includegraphics[scale=0.2]{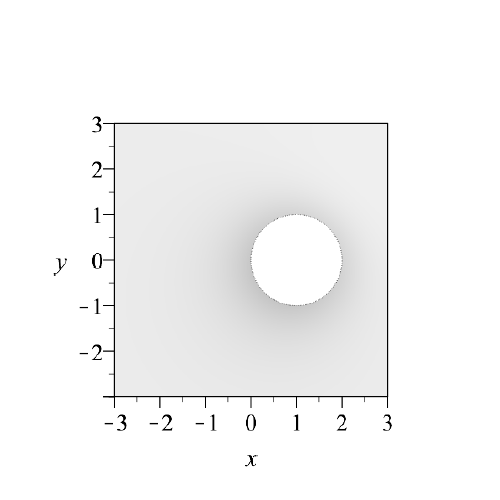}
		\caption{Implicit Euler method}
		\label{fig:theta-methods:c}
	\end{subfigure}
	\caption{The classical stability regions for the one-stage $\theta$ methods (a) $\theta=0$, (b) $\theta=\tfrac{1}{2}$, (c) $\theta=1$. The stability region of (b) is the same for any exactly A-stable method.}
	\label{fig:theta-methods}
\end{figure}

For Runge-Kutta methods, the absolute stability region may be characterized in terms of the quantity $R(\mu)$ as the region satisfying $|R(\mu)|\leq 1$. Having this condition satisfied strictly ($<$) for those problems where the exact solution decays (eigenvalues with negative real part) underlies the criterion of \emph{A-stability}. An A-stable method is any method that contains the right-half plane in its stability region. This is an important criterion for stiff problems, in the sense of problems that have eigenvalues with real parts widely separated and negative, because an A-stable method will avoid stability restrictions that force an extremely small time step when small time steps are not required for an accurate solution. An exactly A-stable method is one that has precisely the left-half plane, with the imaginary axis as boundary, as its absolute stability region. The stability regions of the Euler method, the implicit Euler method and the implicit midpoint rule are shown in figure~\ref{fig:theta-methods}. The implicit midpoint rule is an example of an exactly A-stable method.

\begin{figure}[th!]
	\centering
	\begin{subfigure}[t]{0.45\textwidth}
		\centering
\includegraphics[scale=0.3]{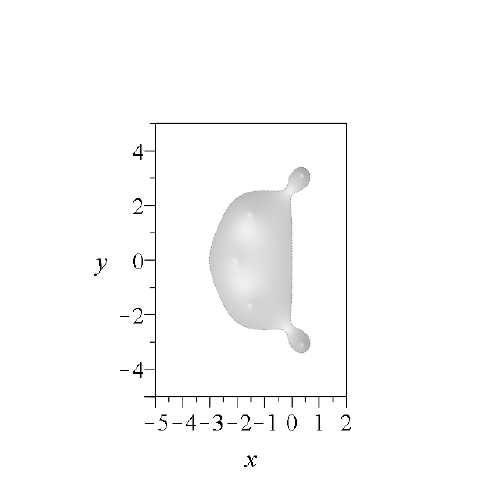}
		\caption{Fourth order component of RKF45}
		\label{fig:rkf45:a}
	\end{subfigure}
	\begin{subfigure}[t]{0.45\textwidth}
		\centering
\includegraphics[scale=0.3]{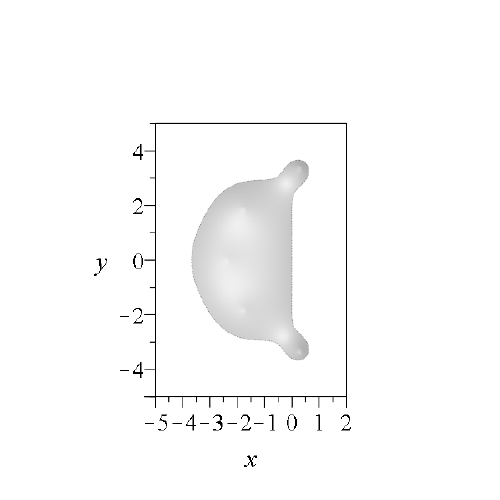}
		\caption{Fifth order component of RKF45}
		\label{fig:rkf45:b}
	\end{subfigure}
	\caption{The classical stability regions for the components of the Fehlberg method.}
	\label{fig:rkf45}
\end{figure}

Higher order methods have more complex and more interesting absolute stability regions. The Runge-Kutta-Fehlberg (or RKF45) method includes a pair of Runge-Kutta methods consisting of a fourth order method together with a fifth order error estimator. The absolute stability regions for these methods are shown in figure \ref{fig:rkf45}.

The shift to consider optimal residual stability regions entails a very different perspective on numerical methods for stiff problems. 

\begin{figure}
\bc
\includegraphics[scale=0.2]{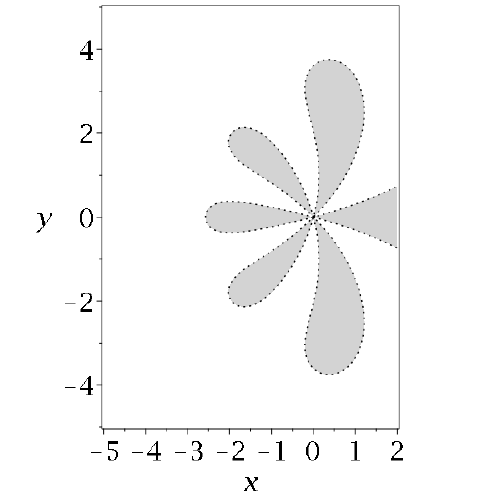}
\ec
\caption{Order star for the fifth order component of RKF45. The $A_{-}$ region is shaded, the $A_0$ boundary is dotted, and the $A_{+}$ region is the unshaded part of the plane.}
\label{fig:order-star}
\end{figure}

As is well-known, consideration of the relative forward error $|R(\mu)e^{-\mu}-1|$ gives rise to the beautiful theory of Order Stars, which considers the regions $A_+: |R(\mu)e^{-\mu}|>1, A_0: |R(\mu)e^{-\mu}|=1, A_-: |R(\mu)e^{-\mu}|<1$. For example, the three regions for the fifth-order method of RKF45 are shown in figure \ref{fig:order-star}. For a detailed account see~\cite{iserles1991order}. In hindsight, it is not a surprise that a similar theory will arise from a consideration of relative backward error.

\begin{table}
\bc
\begin{tabular}{|c|c|}
\hline
method & $R(\mu)$ \\
\hline\hline
Euler & $1+\mu$ \\
\hline
Backward Euler & $\frac{1}{1-\mu}$ \\
\hline
Implicit Midpoint Rule & $\frac{1+\mu}{1-\mu}$\\
\hline
Taylor Series order $p$ & $\sum_{k=0}^p\frac{\mu^k}{k!}$ \\
\hline
\end{tabular}
\ec
\caption{Approximations to $e^{\mu}$ for some numerical methods}
\label{tab:exp-approx}
\end{table}

\subsection{Elementary analysis}
We now consider using only elementary tools to achieve the same result for the Dahlquist test problem, with a one-step method that produces $y_{n+1} = R(\mu)y_n$.
Also without loss of generality consider only $t_n=0$ and $t_{n+1}=h$. We search for an optimal interpolant, which we call $z(t)$, satisfying $z(0)=y_n$, $z(h)=y_{n+1}$ and
\begin{equation}
\label{local-sol}
\dot{z}(t)=\lambda(1+\delta(t))z(t),
\end{equation}
with $\delta(t)$ as small as possible. By constructing it explicitly we will demonstrate its existence and its minimality for $\delta(t)$. We use the $\infty$-norm to measure the size of $\delta$:
$$\|\delta(t)\|_{\infty}=\max_{0\leq t\leq h}|\delta(t)|.$$
Rearranging, and assuming $z(s)\neq 0$ anywhere in $0\leq s\leq h$,
$$\frac{\dot{z}(t)}{z(t)}-\lambda=\lambda\delta(t),$$
so
$$\int_0^{h}\frac{\dot{z}(s)}{z(s)}ds - \lambda\int_0^{h}ds=\lambda\int_0^{h}\delta(s)ds,$$
or
\begin{equation}
\ln z(h)-\ln z(0) - \lambda h = \lambda\int_0^{h}\delta(s)\,ds -2\pi ik\>,
\label{eq:firstunwind}
\end{equation}
where the integer $k$ is determined by the (as yet unknown) number of times the path $z(s)$ winds around $z=0$. We will see later that $k=0$ often, particularly for small step sizes, but there are important cases where $k\neq 0$, especially for higher-order methods. By adjusting $k$ if necessary we get
$$\ln_k\left(\frac{z(h)}{z(0)}\right) - \mu = \lambda\int_0^{h}\delta(s)\,ds,$$
where we have used David Jeffrey's compact notation $\ln_kz$ for $\ln z+2\pi ik$.\\

Since $\frac{z(h)}{z(0)}=\frac{y_{n+1}}{y_n}=R(\mu)$ we have
$$\ln_k R(\mu) - \mu = \lambda\int_0^{h}\delta(s)ds.$$
Taking absolute values and using the triangle inequality,
$$\left|\ln_k R(\mu) - \mu\right| \leq |\lambda|\,\left|\int_0^{h}ds\right|\,\|\delta(t)\|_{\infty}$$
or
$$|\ln_k R(\mu)-\mu|\leq |\mu|\,\|\delta\|_{\infty}$$
or
$$\|\delta\|_{\infty}\geq \left|\frac{\ln_k R(\mu)}{\mu}-1\right|,$$
where we ignore the uninteresting case of $\mu=0$.

This fundamental inequality gives a lower bound on any backward error $\delta(t)$ capable of taking $y_n$ to $y_{n+1}$. We now show that this lower bound is achieved, if we choose $\delta(t)$ to be constant. Suppose
\begin{equation}
\delta=\frac{\ln_kR(\mu)}{\mu}-1\>. \label{eq:optimalresidual}
\end{equation}
Then $\dot{z}=\lambda(1+\delta)z$ is well-defined and indeed
$$z=e^{\lambda(1+\delta)t}z(0)=e^{\lambda\frac{\ln_kR(\mu)}{\mu}t}y_n.$$
When $t=h$,
\bea
z(h)&=&e^{\lambda\frac{\ln_kR(\mu)}{\mu}h}y_n\nonumber\\
&=&e^{\mu\frac{\ln_kR(\mu)}{\mu}}y_n\nonumber\\
&=&e^{\ln_kR(\mu)}y_n\nonumber\\
&=&R(\mu)y_n\nonumber\\
&=&y_{n+1}
\eea
as desired, for any choice of $k$.\\

Since the set $\left\{\left|\frac{\ln_kR(\mu)}{\mu}-1\right|\, :\, k\in\mathbb{Z}\right\}$ is countable and nonnegative, it has a least member. We will use the $k$ that picks out the least member of this set.

\noindent
\emph{Lemma:} To minimize $\left|\frac{\ln_kR(\mu)}{\mu}-1\right|$ over choices of branch cut, we must choose $$k=\overline{\mathcal{K}}:=\left[\frac{\text{Im}(\mu-\ln R(\mu))}{2\pi}\right],$$
where $[a]$ is the nearest integer to $a$; in case of a tie, either integer above or below will do.

\noindent
\begin{proof}
We seek a choice of $k$ that minimizes $\left|\frac{\ln_kR(\mu)}{\mu}-1\right|$, \ie, that minimizes $|\varepsilon|$ for
$$\varepsilon:=\frac{\ln_kR(\mu)}{\mu}-1=\frac{\ln R(\mu)+2\pi ik}{\mu}-1.$$
Putting this over a common denominator, we get
\begin{equation}
    \varepsilon = \frac{ \ln R(\mu) - \mu + 2\pi i k }{\mu}\>.\label{eq:unwindkround}
\end{equation}
For the magnitude of $\varepsilon$ to be as small
as possible, we choose $k$.  But this alters only the imaginary part of the numerator, because
$k$ is an integer and therefore $2\pi ik$ is purely imaginary.
In order to make this imaginary term as small as possible, it follows that we must have that the integer $k$ cancels as much as possible of the imaginary part
from equation~\eqref{eq:unwindkround}, or
\begin{equation}
k=\overline{\mathcal{K}}:=\left[\frac{\text{Im}(\mu-\ln R(\mu))}{2\pi}\right]\>. \label{eq:minimalk}    
\end{equation}

This formula is very reminiscent of the unwinding number from~\cite{jeffrey1996unwinding}, but is different in detail.
\end{proof}
\noindent
\begin{remark}
We will be able to examine the optimal backward error of a very large class of methods, essentially all one-step methods, with just formula~\eqref{eq:optimalresidual} with the integer $k$ specified by formula~\eqref{eq:minimalk}.
\end{remark}

\noindent
\begin{remark} If ever an interpolant hits $z(s)=0$, from $\dot{z}=\lambda(1+\delta)z$ we would have $z\equiv 0$ from then on, so unless $y_{n+1}=0$ we could not reach it. Moreover, if $y_{n+1}=0$ then $\|\delta\|_{\infty}=\infty$ is necessary (a small absolute $\Delta(t)$ is possible in that case, but the relative error must go to infinity). As we saw in the general optimal control approach, this corresponds to the case when the control is singular.
\end{remark}

\s{Results: Optimal backward error for some methods}
We now examine properties of the optimal backward error along with contour plots of the optimal $\infty$-norm optimal residuals. Note $|\delta|=1$ corresponds to $100$\% error; $|\delta|$ larger than this means we are solving a totally different equation. In the absence of systematic structure preservation, we may regard any likeness of the solution to what we want as coincidence.\\ 

Taking $|\delta|=0.05$ is analogous to the ``$95\%$ confidence limit''. If $|\delta|\leq 0.05$, we have the exact solution to a problem within $5$\% of the one we wanted to solve. We will also be concerned with the asymptotic limit $|\delta|=\varepsilon$ as $\varepsilon\rightarrow 0$.


\begin{figure}[th!]
	\centering
	\begin{subfigure}[b]{0.3\textwidth}
		\centering
\includegraphics[scale=0.2]{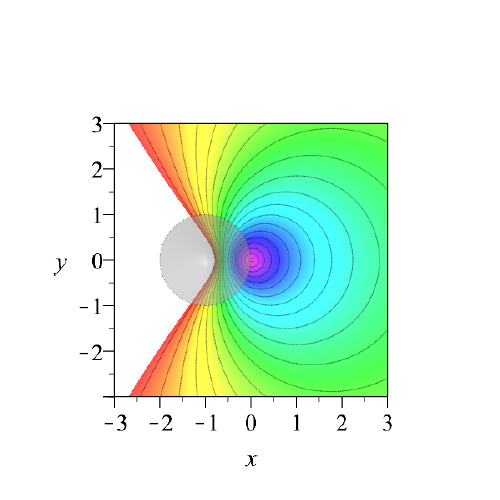}
		\caption{Euler method}
		\label{fig:theta-methods-residual:a}
	\end{subfigure}
	\begin{subfigure}[b]{0.3\textwidth}
		\centering
\includegraphics[scale=0.2]{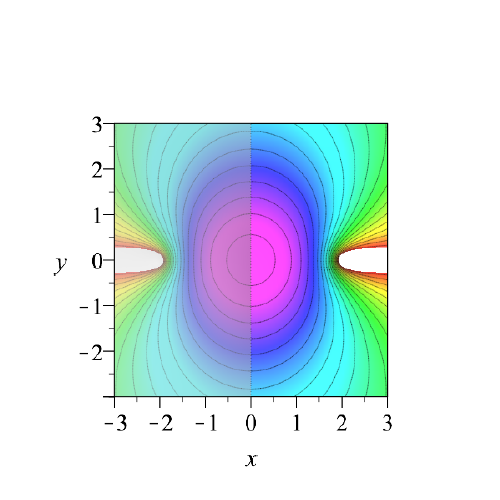}
		\caption{Implicit midpoint rule}
		\label{fig:theta-methods-residual:b}
	\end{subfigure}
	\begin{subfigure}[b]{0.3\textwidth}
		\centering
\includegraphics[scale=0.2]{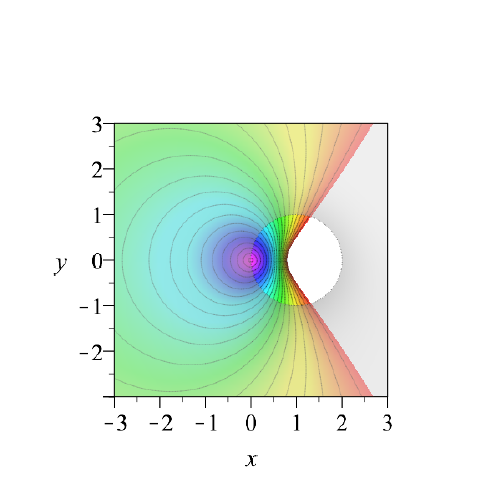}
		\caption{Implicit Euler method}
		\label{fig:theta-methods-residual:c}
	\end{subfigure}
	\caption{The $\infty$-norm optimal residual stability regions from equation~\eqref{eq:optimalresidual} for the one-stage $\theta$ methods (a) $\theta=0$, (b) $\theta=\tfrac{1}{2}$, (c) $\theta=1$. In each image, the classical stability region is shaded and the $\infty$-norm optimal residual stability region is coloured, with contours at 5\% intervals. Regions where $|\delta|>1$ are coloured white. Note the significantly larger
region where $\delta=0.05$ for the second-order method with $\theta=\tfrac12$.  In essence this is why second-order methods are more efficient than first-order methods, for a given accuracy. Perhaps most important, the regions where $|\delta|>1$, \emph{i.e.} those regions where the problem solved is more than $100\%$ different to the problem we wanted to solve, has a nontrivial intersection with the classical stability region.}
	\label{fig:theta-methods-residual}
\end{figure}

\subsection{Theta methods}
The rational approximation that arises from the theta-method
\begin{equation}
y_{n+1} = y_n + \mu\left( (1-\theta) y_n + \theta y_{n+1} \right)
\end{equation}
is
\begin{equation}
R(\mu) = \frac{1+(1-\theta)\mu}{1-\theta \mu}\>.
\end{equation}
For $\theta = 0$ we have explicit Euler, for $\theta = 1$ we have implicit Euler,
and for $\theta = 1/2$ we have the implicit midpoint rule.  The contours of $\delta$
for each of these three $\theta$ are plotted in figure~\ref{fig:theta-methods-residual}.  We see that for explicit Euler there is a substantial portion of the classical stability region $|\mu + 1|<1$ where the computed solution will have $|\delta| > 1$.  This means that even though the solution will decay, it must be the solution of a problem that is more than $100\%$ different to the original problem.  Likewise, a substantial portion of the right-half plane will have $|\delta|>1$ for implicit Euler; this is
why such numerical methods can suppress (actual) chaos~\cite{corless1991numerical}.

In contrast, the implicit midpoint rule (with $\theta = 1/2$) doesn't have such large $100\%$ error zones (although it does, near the portions of the real axis where $|\mu| >2$).  Instead it has an interestingly larger region where $|\delta| = 0.05$ than either explicit or implicit Euler does.  This is because the method is second-order, and the substantially larger region ``of $95\%$ confidence,'' if you like, is a reflection of the value of a second-order method.  Higher order methods really can be more efficient, and this picture shows why.

\begin{figure}[th!]
	\centering
	\begin{subfigure}[t]{0.45\textwidth}
		\centering
\includegraphics[scale=0.35]{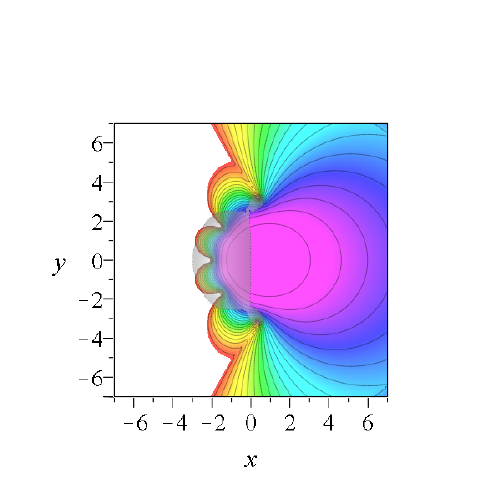}
		\caption{Fourth order component of RKF45}
		\label{fig:rkf45-residual:a}
	\end{subfigure}
	\begin{subfigure}[t]{0.45\textwidth}
		\centering
\includegraphics[scale=0.35]{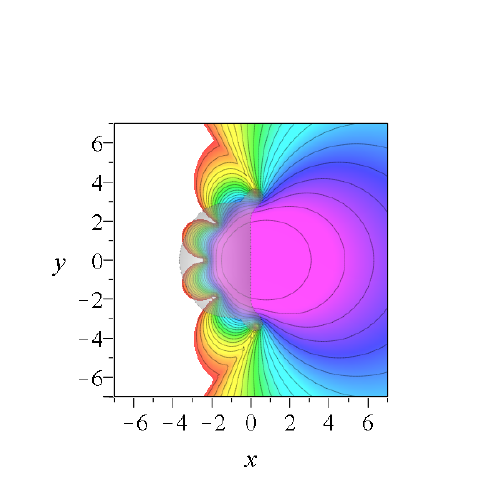}
		\caption{Fifth order component of RKF45}
		\label{fig:rkf45-resiudal:b}
	\end{subfigure}
	\caption{The $\infty$-norm optimal residual stability regions  from equation~\eqref{eq:optimalresidual} for the components of the RK-Fehlberg method. In each image, the classical stability region is shaded and the $\infty$-norm optimal residual stability region is coloured, with contours at 5\% intervals.}
	\label{fig:rkf45-residual}
\end{figure}

\subsection{RKF45}
In figure~\ref{fig:rkf45-residual} we find the optimal
stability regions for each member of the RKF45 pair.  We see substantial regions in the left half plane where the optimal residual must be larger than $1$ in magnitude; that is, where the backward error must be larger than $100\%$.
The regions where $|\delta|\le 0.05$ are comparable for each member of the pair, showing that if one method is accurate, the other is likely to be as well.  This information is complementary to that of the classical stability regions for the pair in figure~\ref{fig:rkf45}.

\begin{figure}[th!]
	\centering
	\begin{subfigure}[t]{0.45\textwidth}
		\centering
\includegraphics[scale=0.35]{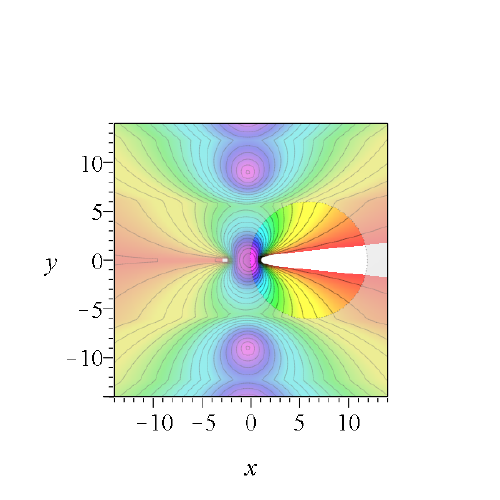}
		\caption{Larger $g$ SDIRK method}
		\label{fig:larger-SDIRK-residual:a}
	\end{subfigure}
	\begin{subfigure}[t]{0.45\textwidth}
		\centering
\includegraphics[scale=0.35]{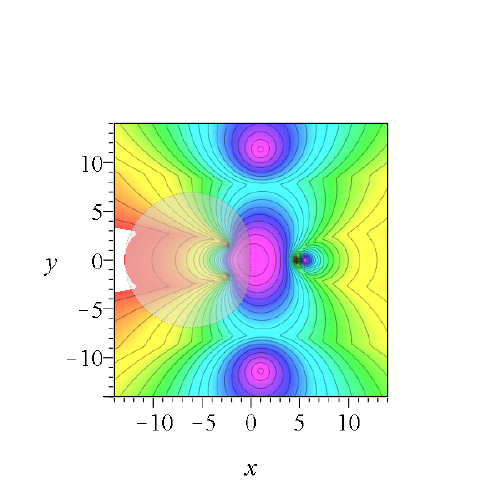}
		\caption{Smaller $g$ SDIRK method}
		\label{fig:smaller-SDIRK-resiudal:b}
	\end{subfigure}
	\caption{The $\infty$-norm optimal residual stability regions  from equation~\eqref{eq:optimalresidual} for the two singly-diagonally-implicit Runge-Kutta (SDIRK) methods of order $3$. In each image, the classical stability region is shaded and the $\infty$-norm optimal residual stability region is coloured, with contours at 5\% intervals. Classical stability analysis suggests that the larger $\gamma$ method (left) is more stable.  This analysis suggests that the smaller $\gamma$ method (right) is more accurate.}
	\label{fig:SDIRK-residual}
\end{figure}

\subsection{A third order SDIRK method}
Singly-Diagonally Implicit Runge-Kutta (SDIRK) 
methods are important and efficient implicit solvers for stiff problems.  In~\cite{wanner1991solving} we find a detailed
derivation of a third-order SDIRK method that
contains the diagonal parameter $\gamma$.  Two values of $\gamma$ make the method $3$rd order. It is argued there that one value of $\gamma$ is to be preferred over the other, as making the method have a larger (much larger) classical stability region.  However, the optimal backward error diagram in figure~\ref{fig:SDIRK-residual} shows that
it is the other value of $\gamma$ that is 
potentially more accurate.

Of course, real life is more complicated than either of these images would show.  See the
massive review~\cite{kennedy2016diagonally} for a more nuanced picture.

\begin{figure}[th!]
	\centering
\includegraphics[width=0.6\textwidth]{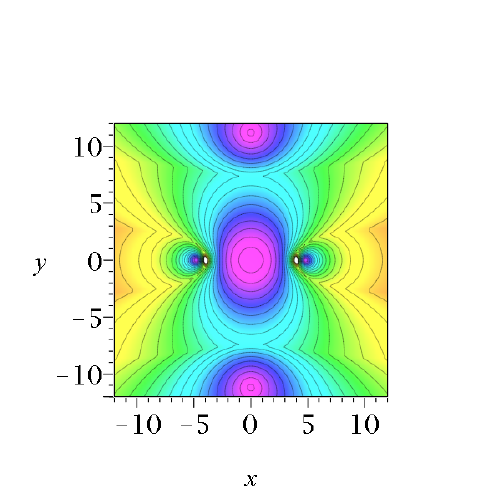}
		\caption{Lanczos $\tau$-method}
	\caption{The $\infty$-norm optimal residual stability regions  from equation~\eqref{eq:optimalresidual} for the Lanczos Tau method with $n=1$. The $\infty$-norm optimal residual stability region is shaded in each image, with contours at 5\% intervals.}
	\label{fig:Tau-residual}
\end{figure}

\subsection{The Lanczos $\tau$-method}
We have included one figure (figure~\ref{fig:Tau-residual}) with the results from solving $\dot y = \lambda y$, $y(0)=y_0$, on $0\le t \le h$, by using the Lanczos $\tau$-method (see for instance~\cite{lanczos1988applied} or~\cite{ortiz1969tau}).  For ease of reference, we include a brief description of the method here.  One starts with a Chebyshev expansion for $\dot y$, namely
\begin{equation}
\dot y(t) = \sum_{k=0}^n c_k T_k(\theta)\>, \label{eq:Lanczosder}
\end{equation}
with $\theta = -1 + 2t/h$.  Integrating $T_0$ with respect to $\theta$ gives $T_1(\theta)$, integrating $T_1(\theta)$ gives $T_2(\theta)/4 + T_0(\theta)/4$,
and thereafter (ignoring constants of integration)
\begin{equation}
\int T_k(\theta)\,d\theta = \frac{1}{2k+2}T_{k+1}(\theta) - \frac{1}{2k-2}T_{k-1}(\theta)\>.
\end{equation}
Lanczos chose to expand the derivative of the unknown,
and then integrate using these simple formulas, because this
was simpler to do by hand.  The details of the process are 
now of course automatable in a variety of ways.
Then integration of equation~\eqref{eq:Lanczosder} gives an expression for $y(t)$, in terms of Chebyshev polynomials of degree $n+1$ and below.
One then uses the initial condition $y(0) = y_0$ to identify the constant introduced on integration.
One subtracts $\lambda$ times this expression from the expression for $\dot y$, and sets the coefficients of $T_k(\theta)$ to zero for $k=0$ up to $k=n$.
This leaves here a residual term containing $T_{n+1}(\theta)$.  If we do this using $n=1$ we find that our Chebyshev interpolant is
\begin{equation}
z(t) = y_0\left( \frac{16}{(4-\mu)^2}T_0(\theta) + \frac{8\mu}{(4-\mu)^2}T_1(\theta) + \frac{\mu^2}{(4-\mu)^2}T_2(\theta)\right)\>,
\end{equation}
which can be simplified at $t=h$ to
\begin{equation}
y_1 = \frac{1 + \mu/2 + \mu^2/16}{1-\mu/2 + \mu^2/16}\,y_0\>.
\end{equation}
This identifies $R(\mu) = (1+\mu/4)^2/(1-\mu/4)^2$ and from there we can identify the optimal $\delta$, which is plotted in figure~\ref{fig:Tau-residual}.
Similar graphs can be produced for larger $n$.

\begin{remark}
The Lanczos $\tau$-method is quite close in spirit to 
finding the minimal $\infty$-norm residual, because 
Chebyshev expansions are, as is well-known, near-minimal for real $\mu$
in this sense. One important difference here is that because
$\lambda \in \mathbb{C}$ we are concerned with the region in the complex plane where $|\delta|$ is small (less than $1$ certainly and less than, say, $0.05$ in practice). Comparison of the optimal backward error on the real interval $(0,h)$ to the size of the Chebyshev residual from the $\tau$ method (not shown here) shows that the Chebyshev method is less than a factor of five worse.  Reluctantly, however, we do not pursue the Lanczos $\tau$-method further here.
\end{remark}

\begin{figure}[th!]
	\centering
	\begin{subfigure}[t]{0.24\textwidth}
		\centering
\includegraphics[scale=0.15]{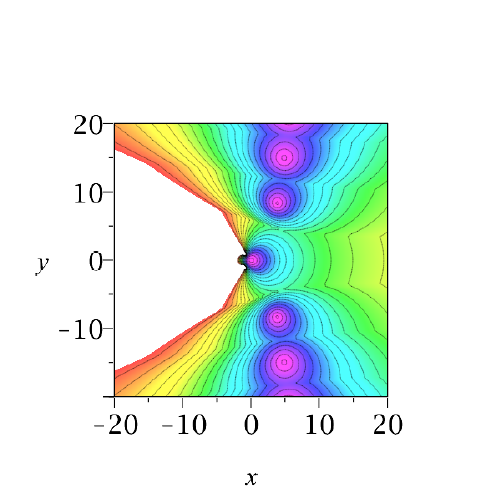}
		\caption{Order 2 Taylor method}
		\label{fig:Taylor-methods-residual:a}
	\end{subfigure}
	\begin{subfigure}[t]{0.24\textwidth}
		\centering
\includegraphics[scale=0.15]{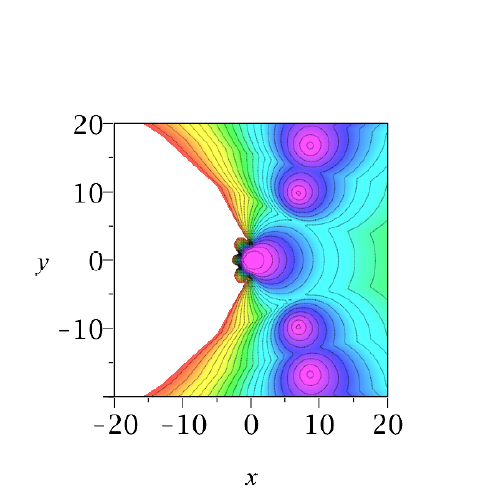}
		\caption{Order 4 Taylor method}
		\label{fig:Taylor-residual:b}
	\end{subfigure}
	\begin{subfigure}[t]{0.24\textwidth}
		\centering
\includegraphics[scale=0.15]{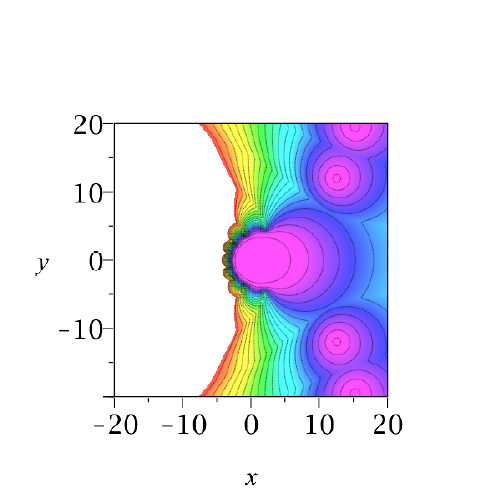}
		\caption{Order 8 Taylor method}
		\label{fig:Taylor-methods-residual:c}
	\end{subfigure}
	\begin{subfigure}[t]{0.24\textwidth}
		\centering
\includegraphics[scale=0.15]{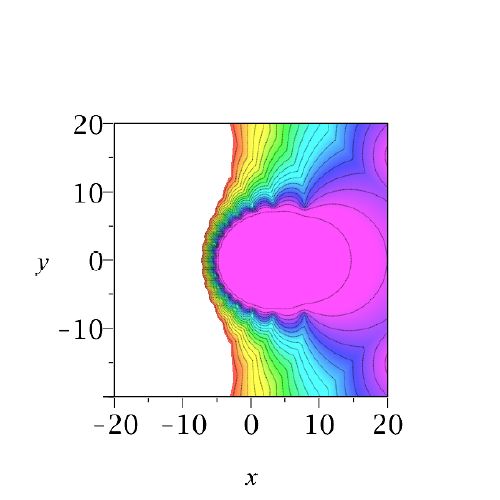}
		\caption{Order 16 Taylor method}
		\label{fig:Taylor-methods-residual:d}
	\end{subfigure}
	\caption{The $\infty$-norm optimal residual stability regions  from equation~\eqref{eq:optimalresidual} for the higher order Taylor methods of order (a) $2$, (b) $4$, (c) $8$ and (d) $16,$. The $\infty$-norm optimal residual stability region is shaded in each image, with contours at 5\% intervals. }
	\label{fig:Taylor-methods-residual}
\end{figure}

\begin{figure}[th!]
	\centering
	\begin{subfigure}[t]{0.24\textwidth}
		\centering
\includegraphics[scale=0.15]{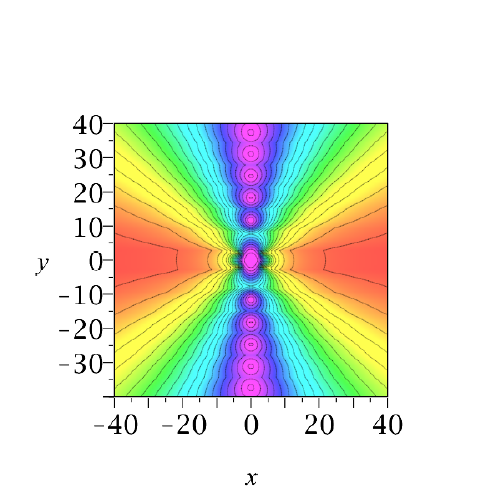}
		\caption{2,2 Pad\'{e} method}
		\label{fig:diag-pade-methods-residual:a}
	\end{subfigure}
	\begin{subfigure}[t]{0.24\textwidth}
		\centering
\includegraphics[scale=0.15]{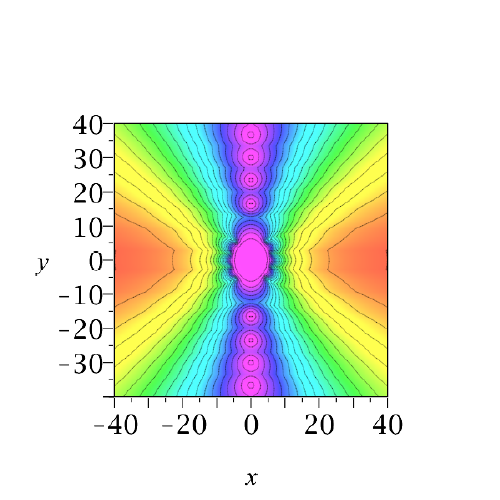}
		\caption{4,4 Pad\'{e} method}
		\label{fig:diag-pade-methods-residual:b}
	\end{subfigure}
	\begin{subfigure}[t]{0.24\textwidth}
		\centering
\includegraphics[scale=0.15]{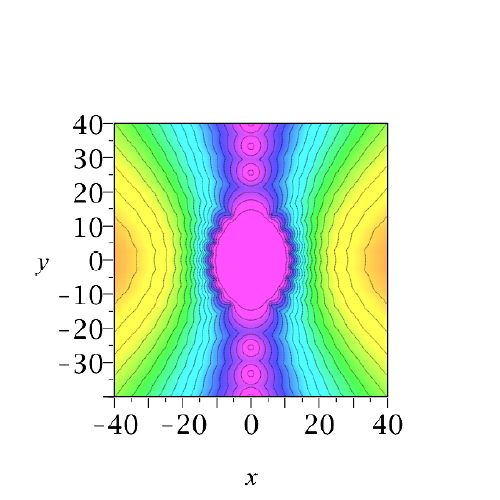}
		\caption{8,8 Pad\'{e} method}
		\label{fig:diag-pade-methods-residual:c}
	\end{subfigure}
	\begin{subfigure}[t]{0.24\textwidth}
		\centering
\includegraphics[scale=0.15]{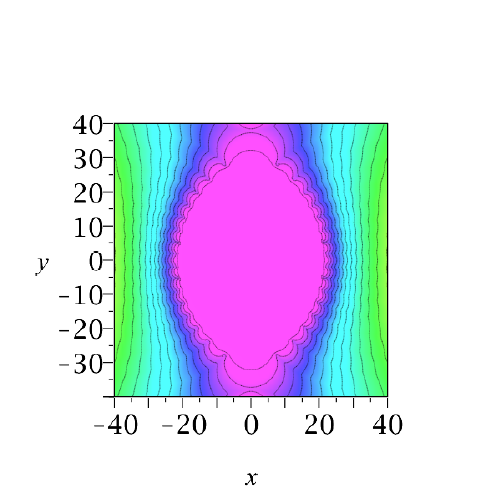}
		\caption{16,16 Pad\'{e} method}
		\label{fig:diag-pade-methods-residual:d}
	\end{subfigure}
	\caption{The $\infty$-norm optimal residual stability regions  from equation~\eqref{eq:optimalresidual} for the higher order Pad\'{e} methods of signature (a) $2,2$, (b) $4,4$, (c) $8,8$ and (d) $16,16$. The $\infty$-norm optimal residual stability region is shaded in each image, with contours at 5\% intervals.}
	\label{fig:diag-pade-methods-residual}
\end{figure}

\subsection{Taylor series methods and Pad\'e methods}
Taylor series methods, including implicit Taylor series
methods, and their generalization to Hermite-Obreshkoff methods, remain of interest for practical problems.
See for instance~\cite{neher2007taylor}.  For the Dahlquist
test problem, all of these methods lead to $R(\mu)$ being
a Pad\'e approximant to $\exp(\mu)$.  Figures~\ref{fig:diag-pade-methods-residual} and~\ref{fig:Taylor-methods-residual} were generated using Maple's built-in function for
computation of Pad\'e approximants~\cite{geddes1979symbolic}.
Notice that as the order of the method increases, the size of the area enclosed by the $\delta = 0.05$ contour increases.  We remark that the unwinding number from formula~\eqref{eq:minimalk} was necessary to get this large central region correct, for high-order methods.  Without the correct branch of logarithm chosen, the central region only had a vertical width of $2\pi$, corresponding to the range of the principal branch of logarithm.

\subsection{Asymptotic results}
Formula~\ref{eq:optimalresidual} is not an asymptotic result.  It does not rely on the time-step $h$ being small.  The optimal $\delta$ is the optimal $\delta$, with almost no reservations or caveats. It 
exists for the Dahlquist test problem for $O(1)$ intervals around $h=0$, regardless of the one-step numerical method used, so long as the optimal interpolant does not go through $0$.  Nonetheless there are some interesting connections to asymptotic results.  First, if the underlying method has forward error $O(h^p)$ as $h \to 0$, that is, it is a $p$th order method, then the residual $\delta$ will also go to zero like $h^p$.  That is, $\delta = O(h^p)$.  This is in contrast to the distinction between so-called local error (which is $O(h^{p+1})$) and forward error (also called global error). For instance, if 
\begin{equation}
R(\mu) = \frac{(1+\mu/4)^2}{(1-\mu/4)^2}\>,
\end{equation}
which it is for the Lanczos $\tau$-method with $n=1$, then 
\begin{equation}
\delta = \frac{\ln_k(R(\mu))}{\mu} - 1 = {\frac{1}{48}}{\mu}^{2}+{\frac{1}{1280}}{\mu}^{4}+O \left( {\mu}^{6}
 \right)
\end{equation}
showing that this particular instance of the Lanczos $\tau$-method acts on the
Dahlquist test problem as if it was a second-order method (in fact it is in general a second order method, but all this expansion does is illustrate that it's second-order on this problem). Note that for small $h$, $R(\mu) = 1 + O(\mu^2)$ and hence the
unwinding number is $k=0$. 

Because the residual is connected to the forward error by the 
Gr\"obner-Alexeev nonlinear variation-of-constants formula~\cite{hairer1993solving}, the forward error is also
$O(h^2)$ for this method.

Similar expansions for other methods confirm their numerical order.
\s{Discussion}
The classical stability analysis using the Dahlquist test problem is informative in that it gives the fundamental reason for loss of stability on taking too large a step size for a stiff problem.  Yet the classical analysis neglects the case of eigenvalues with positive real part (which are indeed important in applications) and furthermore says little about accuracy: decay is a qualitative feature, not a quantitative one, and one can have decay and yet be $100\%$ wrong.  As we see in the figures for some common methods, computing the optimal backward error gives some complementary information, namely the  relative size of the figure enclosed by the $0.05$ contour level: the larger the area, the larger the time-step that can be taken.  We see indeed that this region (and indeed even the $100\%$ region) can be much smaller than the classical stability region.  This indicates that even though the method may be stable for time steps inside the classical stability region but outside the $100\%$ contour, it certainly is not accurate, solving a problem that is more than $100\%$ different to the one that was intended to be solved.

This observation is in some sense not new.  Experienced numerical analysts knew this, and knew that one had to take small time steps for accuracy anyway (when it was needed).  Still, this quantitative assessment of optimal relative backward accuracy is indeed new.  It seems possible that one
might choose a different method to solve some problems, using this criterion, than one might choose using the classical stability criterion.

The analysis here applies to any one-step method. 

\begin{remark}
Since the optimal $\delta$ is constant (for any one-step method), we have the curious observation that constant-stepsize solution using a one-step method 
gives the exact solution of $\dot y = \Lambda y$ where $\Lambda = \lambda(1+\delta) = \lambda(\ln_kR(\mu)/\mu)$; that is, the optimal interpolant is the exact solution of
not only a nearby problem, but a nearby problem of the same kind.  That is, we have not only an optimal backward error, but an optimal \emph{structured} backward error.  

That constant-stepsize one-step methods solved $\dot y = \Lambda y$ was known; what was not known was that this solution has 
the optimality property derived in this paper.
\end{remark}
\begin{remark}
If the optimal backward error is \emph{large}, then the numerical method has necessarily solved a very different problem to the one intended.  This is a sure and certain indication that the underlying numerical method, that generated the skeleton, has failed. In particular, when a numerical method introduces spurious chaos into a nonchaotic system, the optimal backward error must be too large.  A more frequent failure detected (for nonchaotic systems) will be when the
automatically-chosen time step sizes are too large; because the analysis of this paper and the more general paper in preparation do not rely on the asymptotic limit as the mean stepsize goes to zero, detection of failure is sure and certain.
\end{remark}

\s{Concluding Remarks}
John C.~Butcher was correct when, critiquing an earlier incarnation of this idea, he said that the classical stability theory does indeed explain many features of
the behaviour of many numerical methods on stiff problems.  It is a highly-successful theory.  The refinement of using Order Stars is also highly successful, in that new proofs of classical order barriers were obtained using it (for instance).  

However, this optimal backward error idea for $y_{n+1} = R(\mu)y_n$ with $\mu = h\lambda$, namely formula~\eqref{eq:optimalresidual} with the unwinding number~$k$ from formula~\eqref{eq:minimalk}, 
\begin{equation}
\delta=\frac{\ln_kR(\mu)}{\mu}-1\>,
\end{equation}
with
\begin{equation}
k=\overline{\mathcal{K}}:=\left[\frac{\text{Im}(\mu-\ln R(\mu))}{2\pi}\right]\>,   
\end{equation}
also explains, quantitatively, some features of the numerical solution of stiff problems by one-step methods.
In particular, it explains clearly why one might sometimes wish to take smaller timesteps than strictly necessary for stability reasons.  That is, any adaptive step-size control that looks at accuracy must be affected by these regions.  We note that asymptotically as $h \to 0$, the optimal residual approaches the local error per unit step, as explained by a theorem of Stetter (this is cited and extended in~\cite{corless2013graduate}). There are other puzzles that this approach might help to explain.  We look forward,
for instance, to finding an example problem for which the smaller $\gamma$ SDIRK method performs better than the larger $\gamma$ SDIRK method,
as is predicted to exist by figure~\ref{fig:SDIRK-residual}.

\goodbreak
\medskip\noindent
\subsection*{Maple Worksheet}
The \textsc{Maple} worksheet written by R.H.C.~Moir used to produce the
graphs in this paper is at
\url{publish.uwo.ca/~rcorless/Dahlquist/ResidualAnalysis-3D-Contours.mw}.

\subsection*{Acknowledgements}
Part of this work was supported by the Natural Sciences and Engineering Research Council of Canada.  Thanks are also due to the Rotman Institute of Philosophy, the Fields Institute for Research in the Mathematical Sciences, the Ontario Research Center for Computer Algebra, and the University of South Australia which supported a visit of the first author to the second.  Special thanks are owed to John C.~Butcher for his hosting the beautiful ANODE conferences over the years, and his long-term support of research in the Runge-Kutta world.
\small
\bibliographystyle{spmpsci}  
\bibliography{bib}
\s*{}
\end{document}